\documentclass[envcountsect]{svjour3}
\smartqed
\usepackage[usenames,dvipsnames]{color}
\usepackage{bm}
\usepackage{verbatim}
\usepackage{pgfplots}
\pgfplotsset{compat=newest}
\pgfplotsset{plot coordinates/math parser=false}
\usepackage{graphicx}
\usepackage{amsmath}
\usepackage{amssymb}
\usepackage{mathtools}
\usepackage{mathrsfs}
\usepackage{hyperref}
\usepackage{nomencl}
\usepackage{algorithm,algorithmic,multicol}
\usepackage{caption}
\usepackage{slashbox}
\usepackage{subcaption}
\usepackage{makecell}

\newcommand{\argmin}{\operatorname{ argmin}}
\newcommand{\argmax}{\operatorname{ argmax}}

\newcommand{\eg}{{\it e.g.}}
\newcommand{\ie}{{\it i.e.}}
\newcommand{\reals}{{\mathbb{R} }}

\newcommand{\Hil}{{\mathcal{H} }}

\newcommand{\N}{{\mathbb{N} }}

\newcommand{\dist}{{\mathbf{dist} }}

\newcommand{\prox}{\operatorname{\bf prox}}
\newcommand{\refl}{\operatorname{\bf refl}}

\makeatletter

\let\emptyset\varnothing

\newtheorem{assumption}{Assumption}

\DeclareMathAlphabet\mathbfcal{OMS}{cmsy}{b}{n}

\begin{document}

\title{An Inertial Parallel and Asynchronous Forward-Backward Iteration for Distributed Convex Optimization}

\subtitle{}

\author{Giorgos Stathopoulos   \and  Colin N. Jones }

\institute{Giorgos Stathopoulos, Corresponding author \at
             Laboratoire d'Automatique \\
              \'Ecole Polytechnique F\'ed\'erale de Lausanne (EPFL)\\
              Lausanne, Switzerland\\
              georgios.stathopoulos@epfl.ch
           \and
           Colin N. Jones    \at
             Laboratoire d'Automatique \\
              \'Ecole Polytechnique F\'ed\'erale de Lausanne (EPFL)\\
              Lausanne, Switzerland\\
              colin.jones@epfl.ch
}

\date{Received: date / Accepted: date}

\maketitle

\begin{abstract}
Two characteristics that make convex decomposition algorithms attractive are simplicity of operations and generation of parallelizable structures.
In principle, these schemes require that all coordinates update at the same time, \ie, they are synchronous by construction.
Introducing asynchronicity in the updates can resolve several issues that appear in the synchronous case, like load imbalances in the computations or
failing communication links. However, and to the best of our knowledge,
there are no instances of asynchronous versions of commonly-known algorithms combined with inertial acceleration techniques.

In this work we propose an inertial asynchronous and parallel fixed-point iteration from which several new versions of existing convex optimization
algorithms emanate. Departing from the norm that the frequency of the coordinates' updates should comply to some prior distribution, we propose a scheme where the only requirement is that the coordinates update within a bounded interval. We prove convergence of the sequence of iterates generated by the scheme at a linear rate. One instance of the proposed scheme is implemented to solve a distributed optimization load sharing problem in a smart grid setting and its superiority with respect to the non-accelerated version is illustrated.
\end{abstract}


\section{Introduction}\label{intro} 

In recent years, there has been a pronounced interest in revisiting the family of convex optimization algorithms commonly known as \emph{operator splitting schemes} or \emph{decomposition methods}.
This reawakened interest can be primarily attributed to the flourishing of machine learning and like fields, where data sets of unprecedented size are processed. As the result of the enormous computational
workload, parallel computing solutions are often assigned to the task.

The major advantage of operator splitting schemes is per iteration operations that are typically cheap, which is mainly the reason that they are preferred for large scale applications. On the downside, however, are their
sublinear convergence rates, achieving initially quick progress towards some optimal point that subsequently levels off. Several ways to remedy this behavior have been proposed, including alternative metric selections~\cite{combettes2014variable}, relaxation strategies and inertial acceleration~\cite{AlvarezAttouch2001}.

In addition, the inherent parallelization potential of these schemes has spurred a significant amount of research in asynchronous implementations. Asynchronous parallel methods have been mostly motivated from memory allocation applications, when, \eg, a vector is stored in the shared memory space of a multicore computer and can be accessed and altered by the cores in an intermittent manner~\cite{liu2015asynchronous,peng2015arock}.

In this work, we focus on another application area that motivates asynchronicity, namely the existence of an inhomogeneous mixture of agents, where their local updates need not occur at a common rate. This type of problem appears in a setting different from the machine learning ones, \ie, in \emph{multi-agent distributed optimization problems, usually at the presence of a global coordinator}. As an example, in a smart grid setting with distributed resources (agents) and a central operator (coordinator), the local update of a particular agent is the solution to an optimization problem of different complexity than other local subproblems of different agents. In addition, the agents' updates need not occur uniformly, or as a matter of fact, need not draw from any stationary distribution since intermittent failures and delays occur. This would require that the computations of different subproblems are initiated at different time instances and that the agents communicate their solutions to the coordinator in arbitrary sequences. Asynchronous schemes like the one described above pose several challenges in terms of proving convergence in comparison to their synchronous counterparts (see~\cite{bertsekas1989parallel,wright2015coordinate} for interesting overviews).

Our work brings together acceleration techniques with asynchronous implementations of a rather wide family of operator splitting schemes, this of \emph{forward-backward splitting methods (FBS)~\cite[Chapter~25]{book_comb}}. More specifically, we devise an asynchronous iteration in which the coordinates update with varying, arbitrary frequencies and, under some common assumptions, we show that \emph{the distance to the set of fixed points of an inertial version of this asynchronous FBS iteration will converge linearly to zero provided that all the coordinates are visited at least once in a given (bounded) time interval}.

The outline of this paper is as follows: In Section~\ref{sec::notation} we briefly introduce the framework within which we perform the subsequent analysis, namely that of monotone operators. Existing works regarding relaxed and/or inertial fixed-point iterations are visited in Section~\ref{sec::related}. The problem of interest is first formulated and explained in Section~\ref{sec::problem}, where the contributions of this work are also outlined. Section~\ref{sec::proof} illustrates a sketch of the convergence proof of the proposed scheme, while the detailed steps are presented in the Appendices. In Section~\ref{sec::algorithms} we draw the connections of our scheme to existing algorithms and how it gives rise to new versions of the latter. Finally, Section~\ref{sec::example} illustrates the performance of the method in comparison to its regular counterpart for a load sharing problem in the context of a smart distribution grid. Complementary proofs to several sections are provided in Appendix~\ref{app::properties}.

\section{Definitions and Notation}\label{sec::notation} 

The mathematical formalism used for analyzing the algorithms throughout this work stems from monotone operator theory. In this section, we give the most basic definitions and properties that are necessary for understanding the material, which can be found in several sources, \eg,~\cite{book_comb,RyuBoydMonotonePrimer}.

A \emph{relation} or \emph{operator} $T$ in $\Hil$ is a subset of $\Hil\times \Hil$, where $\Hil$ is a Hilbert space. We write $u=T(x)=Tx$ to denote the set $\left\{u\;|\;(x,u)\in T\right\}$. The relation $T$ has \emph{Lipschitz constant $L$} if for all $u\in T(x)$ and $v\in T(y)$ it holds that $\|u-v\|\le L\|x-y\|$. If $L<1$ we call $T$ a contraction, while if $L\le 1$, then $T$ is called \emph{nonexpansive}. A point $x$ is a \emph{fixed point} for $T$ if $x=Tx$, denoted as $x\in\operatorname{fix} T$. This is equivalent to $x\in\operatorname{zer} S$, where $S=I-T$ and $I$ is the identity operator such that $I=\left\{(x,x)\;|\;x\in\Hil\right\}$. The set of fixed points of a nonexpansive operator with full domain ($\operatorname{dom} T=\Hil$) is closed and convex. The operator is called \emph{averaged} if $T=(1-\theta)I+\theta G$ with $\theta\in(0,1)$, where $G$ is a nonexpansive operator. It follows that $T$ is nonexpansive and has the same fixed points as $G$. Compositions of nonexpansive operators are nonexpansive, as well as compositions of averaged operators result in averaged operators.

Monotonicity is another important property of a relation. $T$ is \emph{monotone} if $\langle x-y,Tx-Ty\rangle\ge 0$ for all $x,y\in\Hil$. The relation is \emph{maximal monotone} if there is no monotone operator that properly contains it.

The following definitions are going to play an important role in the rest of the paper.
\begin{definition}{Cocoercivity.}
The operator $T$ is $\delta$-cocoercive with $\delta>0$ if
\[
\langle x-y,Tx-Ty\rangle\geq \delta\|Tx-Ty\|^2,\;\forall x,y\in\Hil.
\]
\end{definition}
\begin{definition}{(Quasi-)Strong monotonicity.}
The operator $T$ is $\mu$-strongly monotone with $\mu>0$ if
\[
\langle x-y,Tx-Ty\rangle\geq \mu\|x-y\|^2,\;\forall x,y\in\Hil.
\]
If the inequality holds only for $y\in\operatorname{zer} T$, \ie,
\[
\langle x-y,Tx\rangle\geq \mu\|x-y\|^2,\;\forall x\in\Hil,
\]
then the operator is quasi-$\mu$-strongly monotone.
\end{definition}

Two operators that play a critical role in finding a zero of a relation are the \emph{resolvent} of $T$, defined as $J_T=(I+T)^{-1}$, and the \emph{reflection} of $T$ defined as $R_T=2J_T-I$. The following properties hold:
\begin{itemize}
 \item If $T$ is monotone, then $J_T$ and $R_T$ are nonexpansive functions.
 \item If $T$ is maximal monotone, then $J_T$ and $R_T$ have full domain $\Hil$.
 \item $T$ shares the same fixed points with $J_T$ and $R_T$.
\end{itemize}

In convex optimization problems, we are interested in specific operators, the most important of which is the \emph{subdifferential $\partial f$} of a convex closed and proper function $f: \Hil\to\reals\cup\{+\infty\}$. The class of convex closed proper functions from $\Hil$ to $\reals\cup\{+\infty\}$ is denoted hereafter with $\Gamma_0(\Hil)$. If $f\in\Gamma_0(\Hil)$, then $\partial f$ is maximal monotone.

The subdifferential is important because the solution to a convex optimization problem can be cast as finding a zero of the subdifferential, \ie, $0\in\partial f(x)$. This can be equivalently written as $x\in(I+\partial f)(x)\Leftrightarrow x\in(I+\partial f)^{-1}(x)\Leftrightarrow x=J_{\partial f}(x)$. Hence the set of optimizers coincides with the fixed points of the resolvent and the reflection of the subdifferential. These operators take special forms in the framework of convex optimization. More specifically, we define the \emph{proximal operator} $\prox{_{\gamma f}}:\Hil\to\Hil$ as the resolvent of the subdifferential of $f$ evaluated at $y$ as $\prox{_{\gamma f}}(y)=(I+\gamma\partial f)^{-1}(y)=\underset{x}{\argmin}\left\{f(x)+(1/2\gamma)\|x-y\|^2\right\}$, $\gamma>0$. Accordingly, the reflection operator is denoted as $\refl{_{\gamma f}}:\Hil\to\Hil$ and is defined as $\refl{_{\gamma f}}=2\prox{_{\gamma f}}-I$.

Throughout this work, a sequence is indexed by a subscript, \ie, $x_{k}$, which corresponds to \emph{the algorithmic iterate $k$}. We use brackets in order to refer to either a single coordinate, or a group of coordinates of the sequence, \ie, $x_k[i]$ refers to the $i^\mathrm{th}$ (block of) coordinate(s) of the sequence $x$ at iteration $k$. When an operator $T$ acts on a point $x$, we might refer to the $i^\mathrm{th}$ coordinate of the output as $(Tx)[i]$.

\section{Related work}\label{sec::related} 

\subsection{Inertial acceleration: The Heavy Ball Method. }\label{ssec::HBM} 
The celebrated Heavy Ball Method (HBM), developed by Polyak in the seminal work~\cite{PolyakHBM}, is a modification of the gradient descent iteration that generates a sequence of iterates $\{x_k\}_{k\in\N}$ that minimize a differentiable, convex function $f$
\begin{equation}{\label{eq::HBM}}
 x_{k+1} = x_k - \gamma\nabla f(x_k) + \beta(x_k-x_{k-1})\enspace,
\end{equation}
with $\gamma$ being an (admissible) stepsize and $\beta\in(0,1)$.
The algorithm is very similar to gradient descent, with the addition of an extrapolation sequence $\{x_k-x_{k-1}\}_{k\in\N}$ that alters the direction in which the new iterate will be landed by injecting some previous information. This term is also called an `inertial' or a `momentum' term.

The method has been significantly improved in order to tackle more general problems, leading to the appearance of projected, proximal~\cite{ipiasco} as well as incremental~\cite{gurbuzbalaban2015convergence} variants. In addition, its convergence rate has been studied and analyzed. More specifically, the method can be shown to converge linearly~\cite{PolyakBook} when the function $f$ is both smooth (Lipschitz continuous gradient) and strongly convex. In the case that the latter assumption is dropped, and under suitable choices for the parameters $\gamma$ and $\beta$, a $\mathcal{O}(1/k)$ ergodic global convergence rate has been proven in~\cite{ghadimi2015global}.

The method was generalized in the context of finding a zero of a maximal monotone operator $S$
\begin{equation*}
 \mbox{find } x_\ast\in\Hil \; \mbox{ such that }\; 0\in S x_\ast\enspace.
\end{equation*}

In~\cite{AlvarezAttouch2001}, the authors proposed an Inertial-Prox algorithm that generalizes~\eqref{eq::HBM} to
\begin{equation}{\label{eq::ProxIntert}}
 x_{k+1} = J_{\gamma_k S}(x_k + \beta_k(x_k-x_{k-1}))\enspace,
\end{equation}
where $J_{\gamma_k S}$ is the resolvent of $S$. Iteration~\eqref{eq::ProxIntert} generalizes the proximal point algorithm and finds a zero of a maximal monotone operator $S$ by making use of the momentum term. In~\cite{Moudafi2003447}, the authors extended the inertial scheme~\eqref{eq::ProxIntert} to find a zero of the sum of two maximal monotone operators.

\subsection{The Krasnosel'ski\u{i}-Mann iteration. }\label{ssec::KM} 
The Krasnosel'ski\u{i}-Mann (KM) iteration is a \emph{fixed-point iteration} of the form
\begin{equation}{\label{eq::KMiter}}
 x_{k+1} = x_k + \eta_k (Tx_k-x_k)\enspace,
\end{equation}
where $T$ is a nonexpansive operator and $\eta_k\in[0,1],\;\sum_{k\in\mathbb{N}}\eta_k(1-\eta_k)=+\infty$ is a relaxation constant. The KM iteration converges to a solution if one exists~\cite{Krasnoselskii,Mann},~\cite[Theorem~5.14]{book_comb}.

Relaxing a fixed-point iteration is yet another way to speed up the convergence.
Many of the popular operator splitting methods can be cast as a KM iteration, including the Forward-Backward Splitting (FBS) algorithm which we analyze in this work, as well as the Alternating Direction Method of Multipliers (ADMM) and the Douglas-Rachford splitting algorithm (DRS). The derivations can be found in, \eg,~\cite{liang2014convergence}.

The combination of both acceleration techniques, \ie, inertia and relaxation has been proposed in~\cite{Alvarez:2004:WCR}, where the zero of a maximal monotone operator $S$ is recovered by means of the iteration
\begin{equation*}
 x_{k+1} = y_k + \eta_k (J_{\gamma_k S}y_k - y_k)\enspace,
\end{equation*}
where $y_k=x_k + \beta_k(x_k-x_{k-1})$. The iteration converges weakly to a fixed point of $J_{\gamma_k S}$ in a Hilbert space setting.

The works~\cite{AlvarezAttouch2001,Moudafi2003447,Alvarez:2004:WCR} are put under a common framework in~\cite{Maingé2008223}, where the author develops convergence theorems for a generic inertial KM-type iteration of the form
\begin{equation*}
 x_{k+1} = y_k + \eta_k (T_ky_k-y_k)\enspace,
\end{equation*}
with $y_k=x_k + \beta_k(x_k-x_{k-1})$. Convergence is proven under different choices for the parameter sequences $\{\eta_k\}_{k\in\N}$, $\{\beta_k\}_{k\in\N}$ as well as the operator sequence $\{T_k\}_{k\in\N}$.
Finally, in the recent work~\cite{hendrickx} the authors employ relaxed and inertial schemes to accelerate the KM iteration by means of algorithms that auto-tune the involved parameters.

\section{Problem description}\label{sec::problem} 

\subsection{Asynchronous updates. }\label{sec::async} 
The proposed setting involves $N$ agents, each one assigned to update one (block of) coordinate(s) of $x$, \ie, $x=(x[1], \ldots, x[N])\in\Hil$, and $x[i]\in\Hil_i$. The agents seek convergence to a fixed point of a nonexpansive operator $T$. One way to achieve this is to perform block-coordinate updates of the KM iteration~\eqref{eq::KMiter}. Such a scheme has been proposed and analyzed in~\cite{peng2015arock}. The iteration reads
\begin{equation}{\label{eq::coordinateKM}}
 x_{k+1}[i] = x_{k}[i] - \eta_k(Sx_\mathrm{read}^i)[i]\enspace,
\end{equation}
where $S=I-T$, hence the set of fixed points of $T$ is the set of zeros of $S$, \ie, $\operatorname{fix} T = \operatorname{zer} S$.

Iteration~\eqref{eq::coordinateKM} assumes \emph{the existence of a global coordinator, associated with a global clock. All the agents update continuously and in parallel, while the global clock updates the subscript $k$ every time that an agent updates}. The variable $x_\mathrm{read}^i$ \emph{represents the state of the vector $x$ as it existed at the coordinator when the agent that is about to update ($i$) requested it (a `read' operation)}. Each update involves the most recent state of $x$, denoted by $x_k$, and the result of the operator $S$ acting on an outdated version $x_\mathrm{read}^i$. The distinction between $x_k$ and $x_\mathrm{read}^i$ is important, since, on the coordinator's level, several components of $x_k$ have possibly been altered since the time instant that $x_\mathrm{read}^i$ was read. Every global clock count $k$ is uniquely associated to an updated $i^\mathrm{th}$ group of coordinates. In this way, only the $i^\mathrm{th}$ block of rows of the operator $S$ contributes to the next update of $x$, and only the $x[i]\in\Hil_{i}$ block is updated.

Our goal is to propose an accelerated version of~\eqref{eq::coordinateKM} in order to achieve better practical performance without increasing the computational complexity of the iteration.

\subsection{An asynchronous inertial forward-backward iteration. }\label{sec::ourAlg} 
We propose an asynchronous inertial KM iteration scheme for finding a zero of $S$. We confine our interest to not just any KM iteration, but we rather assume
that the operator $T$ can be written as the composition of two operators $T_A:\Hil\mapsto\Hil$ and $T_B:\Hil\mapsto\Hil$, the properties
of which will be analyzed in the course of this section. In addition, we assume that the operator $T_A$ is separable
into $N$ components, \ie, $T_A=(T_{A_1},\ldots,T_{A_N})$, $T_{A_i}:\Hil_i\mapsto\Hil_i$.

The scheme comprises $N+1$ main blocks, one associated to the coordinator and $N$ associated to the agents. The operators $T_{A_i}$ are private to the agents, while $T_B$ is owned by the coordinator. Before proceeding to the algorithm, we introduce the quantities associated to the coordinator and the agents. Let us start by denoting all variables stored at the coordinator with `$x$' and all variables stored at the agents by `$y$'.

\begin{itemize}
\item Coordinator
  \begin{itemize}
  \item $x$ - the current value of the (global) optimization variable
  \item $x^i_\mathrm{write}$ - the value of $x$ at the time of receipt of value from agent $i$
  \item $x^i_\mathrm{read}$ - the value of $x$ at the time of transmission to agent $i$
  \item $z^i$ - the last value received from agent $i$ ($z^i=T_{A_i}\left(y_B+\beta(y_\mathrm{write}-y_\mathrm{write}^\mathrm{prev})\right)$)
  \end{itemize}
\item The following variables are local to agent $i$
  \begin{itemize}
  \item $y_\mathrm{write}=x^i_\mathrm{write}[i]$ - the value of $x[i]$ after updating $x$ with the latest $z^i$
  \item $y_\mathrm{write}^\mathrm{prev}$ - the value of $x[i]$ before $y_\mathrm{write}$
  \item $y_B=(T_Bx^i_\mathrm{read})[i]$ - quantity computed by the coordinator and transmitted at the same time as $y_\mathrm{write}$
  \end{itemize}
\end{itemize}
Agent $i$ essentially waits to receive the quantities $y_B$ and $y_\mathrm{write}$ from the coordinator. Once received, the privately owned operator $T_{A_i}$ is applied to the expression $y_B+\beta(y_\mathrm{write}-y_\mathrm{write}^\mathrm{prev})$ and the result, \ie, $z^i$, is transmitted back to the coordinator, which, in turn, uses it in order to update the $i^\mathrm{th}$ component of the global variable $x$.

The algorithmic scheme can be described by two interacting and distinct blocks, one refering to an agent and one to the coordinator.
\begin{algorithm}[H]
\caption{Agent}
\label{al::AsInFBS_agent}
\begin{algorithmic}
 \STATE \textbf{wait until}
 \STATE \quad Receive $y_\mathrm{write}, y_B$ from coordinator
 \STATE \quad \textbf{Compute:} \begin{equation}{\label{eq::local_solve}}z^i = T_{A_{i}} \left( y_B+\beta(y_\mathrm{write}-y_\mathrm{write}^\mathrm{prev}) \right)\end{equation}
 \STATE \quad Transmit $z^i$ to coordinator
 \STATE \quad $y_\mathrm{write}^\mathrm{prev}\leftarrow y_\mathrm{write}$
 \STATE \textbf{end}
\end{algorithmic}
\end{algorithm}

\begin{algorithm}[H]
\caption{Coordinator}
\label{al::AsInFBS_coordinator}
  \begin{multicols}{3}
    \begin{algorithmic}
      \scriptsize
      \STATE \underline{\textbf{Write Thread}}
      \STATE Initialize $\mathcal{W}=\emptyset$, $\mathcal{R}=\emptyset$, $k=0$
      \REPEAT{}
       \STATE Receive $z^i$ from agent $i$
       \STATE $\mathcal{W}\leftarrow\mathcal{W}\cup\{i\}$
      \UNTIL{stopping condition holds}
    \end{algorithmic}
    \columnbreak
    \begin{algorithmic}
      \scriptsize
      \STATE \underline{\textbf{Compute Thread}}
      \REPEAT{}
      \STATE Choose $i\in\mathcal{W}$ or \textbf{block until} $\mathcal{W}\neq\emptyset$
      \STATE $\mathcal{W}\leftarrow\mathcal{W}\setminus\{i\}$
      \STATE $z[i]\leftarrow z^i$
      \STATE \textbf{Compute:} \begin{equation}{\label{eq::FP_iter}}x_{k+1} = (1-\eta)x_{k} + \eta z\end{equation}
      \STATE $x^i_\mathrm{write}\leftarrow x_k$
      \STATE $\mathcal{R}\leftarrow\mathcal{R}\cup\{i\}$
      \STATE $k\leftarrow k+1$
      \UNTIL{stopping condition holds}
    \end{algorithmic}
    \columnbreak
    \begin{algorithmic}
      \scriptsize
      \STATE \underline{\textbf{Read Thread}}
      \REPEAT{}
      \STATE Choose $i\in\mathcal{R}$ or \textbf{block until} $\mathcal{R}\neq\emptyset$
      \STATE $\mathcal{R}\leftarrow\mathcal{R}\setminus\{i\}$
      \STATE $x^i_\mathrm{read}\leftarrow x_k$
      \STATE \begin{align*}
             &\text{Transmit}  && x^i_\mathrm{write}[i],\\
             &                 && (T_B x^i_\mathrm{read})[i]
             \end{align*}
      \UNTIL{stopping condition holds}
    \end{algorithmic}
  \end{multicols}
\end{algorithm}

\begin{figure}[!htb]
  \begin{center}
  \includegraphics[scale=0.5]{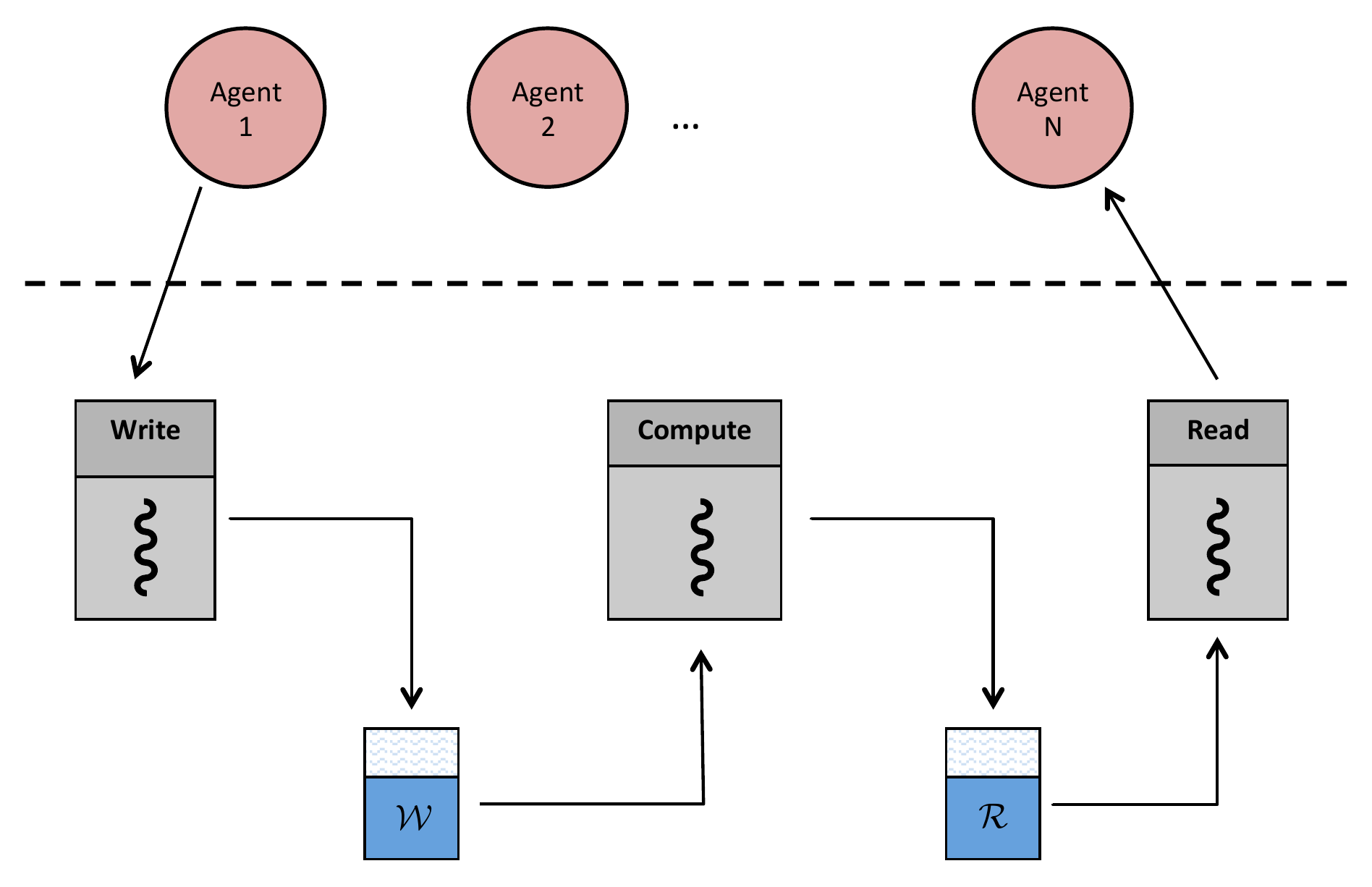}
  \caption{Buffer $\mathcal{W}$ is filled by the Write thread and emptied by the Compute thread. Similarly, $\mathcal{R}$ is filled by the Compute thread and emptied by the Read thread. The scheme executes continuously and asynchronously, both at the coordinator level (concurrent threads) and at the interface between the agents and the coordinator.}
  \label{fig::communication}
  \end{center}
\end{figure}

We make the following observations:
\begin{itemize}
\item Algorithms~\ref{al::AsInFBS_agent} and~\ref{al::AsInFBS_coordinator} are executed continuously and in parallel. There is one algorithmic block described by Algorithm~\ref{al::AsInFBS_agent} per agent $i,\;i=1,\ldots,N$, and these $N$ blocks execute in parallel. Each block is activated \emph{upon receival} of the required info from the coordinator.
\item There are three threads of control in Algorithm~\ref{al::AsInFBS_coordinator}, namely a \emph{Write thread}, a \emph{Compute thread} and a \emph{Read thread}. The threads are \emph{concurrent} and their execution is determined by two buffers, the \emph{read} buffer denoted by $\mathcal{R}$ and the \emph{write} buffer denoted by $\mathcal{W}$.
\item Whenever the coordinator receives an update from an agent, $\mathcal{W}$ is updated. Receivals contribute, therefore, in \emph{filling-up} the write buffer. The coordinator eventually decides to pull an agent from the buffer and use its corresponding value to update $x_{k+1}[i]$ with $z[i]$~\ref{eq::FP_iter}, while the rest of the coordinates are updated based on the previous values of $z[j],\;j\neq i$. Once the update has occured, the index of the corresponding agent is removed from $\mathcal{W}$ and added to $\mathcal{R}$, signaling that the agent is ready to `listen' from the coordinator. Similarly, whenever the coordinator decides to transmit to an agent (Read thread), its index is removed from $\mathcal{R}$. In this way, the buffers control and execution of Algorithm~\ref{al::AsInFBS_coordinator}, which follows a \emph{producer-consumer pattern}.
\item Note that whenever $\mathcal{W}$ is emptied, the Compute thread is blocked until at least one index is added to its stack. The same holds for $\mathcal{R}$ and the Read thread.
\item As indicated by Assumption~\ref{Ass::bouded_delay} below, each agent has to contribute in updating~\eqref{eq::FP_iter} at least once every $\tau$ time epochs. Consequently, the buffers cannot remain empty for longer than $\tau$.
\item The relaxation parameter $\eta$ and the inertia constant $\beta$ that appear in~\eqref{eq::FP_iter} and~\eqref{eq::local_solve}, respectively, will be restricted within intervals in the subsequent sections so as to ensure convergence of the algorithm.
\item The workhorse behind the conceptual scheme derived above is essentialy the asynchronous fixed-point iteration~\eqref{eq::FP_iter} that is based on the agents' updates~\eqref{eq::local_solve}.
\item The coordinator's variables $x,x^i_\mathrm{write},x^i_\mathrm{read}$ and $z$ are containers that get updated from the agents and lie in $\mathcal{H}$. The local variables $y_\mathrm{write}, y_\mathrm{write}^\mathrm{prev}, y_B$ and $z^i$ lie in $\mathcal{H}_i$.
\end{itemize}
Figure~\ref{fig::communication} demonstrates graphically the information flow.

We make the following standing assumptions:
\begin{assumption}{\label{Ass::propertiesA}}
The operator $T_A$ is nonexpansive.
\end{assumption}
\begin{assumption}{\label{Ass::propertiesB}}
The operator $T_B=I-\gamma B$, $\gamma>0$, and the operator $B$ is $1/L$-cocoercive.
\end{assumption}
\begin{corollary}\label{lem::Scocoercive}
The operator $S=I-T$ is $1/2$-cocoercive.
\end{corollary}
The proof can be found in Appendix~\ref{app::Scocoercive}.
\begin{assumption}{\label{Ass::S}}
The operator $S=I-T$ is quasi-$\nu$-strongly monotone for some $\nu>0$.
\end{assumption}
Iteration~\eqref{eq::FP_iter} generalizes the inertial proximal iteration in~\cite{ipiasco}, with $T_A$ replacing the proximal operator and $T_B$ the operator $I-\gamma\nabla f$.
Assumption~\ref{Ass::S} can be met for a relatively wide class of operators $T_A$ and $T_B$. One such instance is derived in Appendix~\ref{app::Scocoercive}, where $\mu$-strong monotonicity of the operator $B$ is assumed in order for the property to hold. In the case of the proximal gradient method with $T_A=\prox_{\gamma g}$ and $T_B=I-\gamma \nabla f$, this assumption would translate to strong convexity of $f$.

Finally, the following assumption regards the frequency of the updates:
\begin{assumption}{\label{Ass::bouded_delay}}
Each agent `writes' to the coordinator state at least once every $\tau$ time epochs.
\end{assumption}
Assumption~\ref{Ass::bouded_delay} categorizes our scheme with the \emph{partially asynchronous parallel methods} as introduced in~\cite[Chapter~7]{bertsekas1989parallel}.

\subsection{Main contribution of the paper. }\label{sec::contribution} 
This paper proves \emph{linear convergence} of the sequence $\{\dist_k\}_{k\in\mathbb{N}}$ generated from Algorithms~\ref{al::AsInFBS_agent} and~\ref{al::AsInFBS_coordinator}, where $\dist_k:=\|x_k-x_\ast\|$. The result is based on Lemma~\ref{lem::ISS} below, which originally appeared in~\cite{lyapunov_approx} and has been extensively used in recent works for proving linear convergence of sequences with errors.
\begin{lemma}{\label{lem::ISS}}
 Let $\{V_k\}_{k\in\N}$ be a sequence of nonnegative real numbers satisfying
 \[
  V_{k+1} \leq rV_k+q\underset{k-\tau\leq l\leq k}{\max}V_l, \quad l\geq 0,
 \]
 for some nonnegative constants $r$ and $q$. If $r+q<1$, then
 \[
  V_k\leq s^kV_0, \quad k\geq 1,
 \]
 where $s=(r+q)^\frac{1}{1+\tau}$.
\end{lemma}
Our convergence proof follows the styles of~\cite{gurbuzbalaban2015convergence} and~\cite{peng2015arock}. In the former, the authors prove linear convergence of the incremental aggregated unconstrained gradient method, while in the latter an asynchronous KM iteration is developed. Our contributions are summarized below.
\begin{enumerate}
 \item \emph{We prove convergence of an asynchronous and parallel forward-backward iteration of a sequence involving an inertial term}. To the best of our knowledge, this is the first result on accelerated asynchronous fixed-point iterations.
 \item Contrary to the majority of existing popular schemes, \emph{the proposed asynchronous iteration is deterministic, \ie, an arbitrary (block of) coordinate(s) can be selected to update at each iteration}. The coordinates can be chosen with varying frequencies, the only assumption being that each coordinate is updated at least once within a fixed time interval. The only work known to us that treats asynchronous updates in a deterministic way, though in a different setting, is~\cite{async_block}.
 \item The proposed iteration is quite general and encompasses many known algorithms as special cases, \ie, several forms of the forward-backward splitting method (see Section~\ref{sec::algorithms}). Indicatively, we propose new inertial and asynchronous instances of two existing and commonly used algorithms and we list some more that can be derived.
\end{enumerate}

\section{Convergence proof}\label{sec::proof} 
We want to use Lemma~\ref{lem::ISS}, with $V_k=\|x_k-x_\ast\|^2,\;x_\ast\in\operatorname{zer} S$. In order to do so, we first express the outdated versions of the global vector $x$ that appear in~\eqref{eq::local_solve} (and consequently in iteration~\eqref{eq::FP_iter}) with respect to the original ones, perturbed by some additive errors. Subsequently, these errors are going to be upper-bounded by $\underset{k-K\leq m\leq k}{\max}\|x_m-x_\ast\|$, for some bounded delay $K$. We are going to go through the proof in steps.

\subsection{Express delayed variables as additive error. }
The variables $y_B,y_\mathrm{write},y_\mathrm{write}^\mathrm{prev}$ that appear in the update~\eqref{eq::local_solve} depend on outdated components of $x$, namely on the state of the vector $x$ when a `read' or `write' operation was performed by agent $i$. It can be easily seen that \emph{any past vector $x_{k-l},\;l\in\{1,\ldots,k-1\}$ can be expressed as}
\[
x_{k-l} = x_k - \sum_{m=k-l}^{k-1}(x_{m+1}-x_m)\enspace.
\]
Consequently, the vectors that appear in~\eqref{eq::local_solve} can be written as functions of the current vector $x_k$ and some error.
\begin{align}{\label{eq::undo_changes_rw}}
x^i_\mathrm{read} \nonumber &= x_{k-l_i} = x_k - a_k^i, \quad \mbox{ for some } l_i\in\{1,\ldots,2\tau\}\\
y_\mathrm{write} \nonumber  &= x^i_\mathrm{write}[i] = x_{k-l_i}[i] = x_k[i] - b_k^i[i], \quad \mbox{ for some } l_i\in\{1,\ldots,2\tau\}\\
y_\mathrm{write}^\mathrm{prev} &= x_{k-l_i}[i] = x_k[i] - c_k^i[i], \quad \mbox{ for some } l_i\in\{1,\ldots,3\tau\}\enspace,
\end{align}
and the sequences $\{a_k^i\}_{k\in\mathbb{N},1\leq i\leq N}$, $\{b_k^i\}_{k\in\mathbb{N},1\leq i\leq N}$ and $\{c_k^i\}_{k\in\mathbb{N},1\leq i\leq N}$ are all of the form $\sum_{m=k-l_i}^{k-1}(x_{m+1}-x_m)$ for some proper choice of $l_i$.
In other words, equations~\eqref{eq::undo_changes_rw} `undo' all the changes that occured over the last updates, until the corresponding past state is recovered. Note that the $l_i$'s in the three equations above are not the same.

\begin{figure}[!htb]
  \begin{center}
  \includegraphics[scale=0.5]{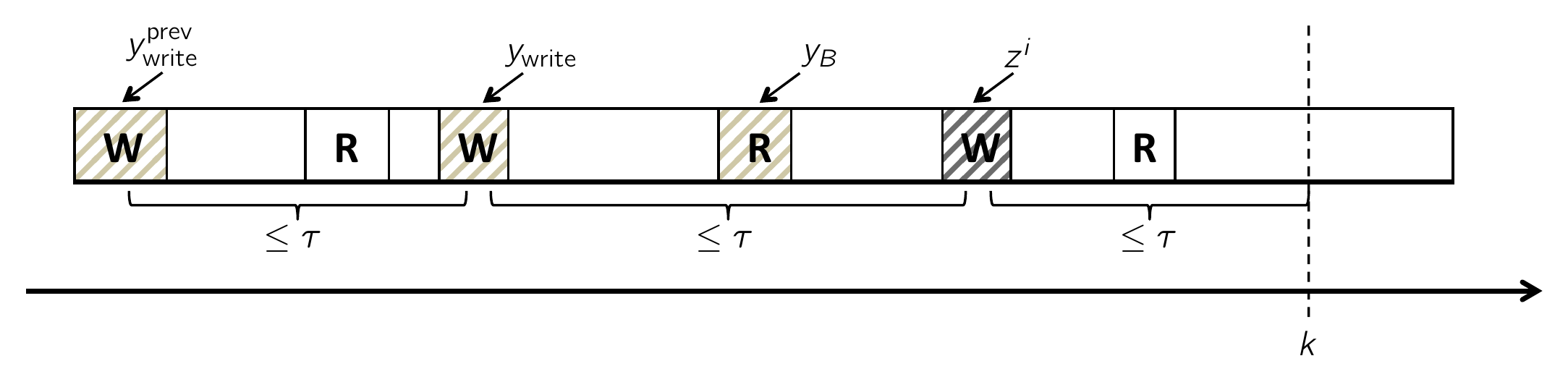}
  \caption{An update is about to occur at $k+1$. From Assumption~\ref{Ass::bouded_delay}, the observed agent will update again no later than $\tau$ time epochs after $z^i$ was communicated to the coordinator. Consequently, $y_\mathrm{write}$ cannot be further than $2\tau$ from the next update, while $y_\mathrm{write}^\mathrm{prev}$ cannot be further than $3\tau$.}
  \label{fig::delays}
  \end{center}
\end{figure}

The intervals within which the subscripts $l_i,\;i=1,\ldots,N$ reside are derived based on Assumption~\ref{Ass::bouded_delay}. The derivation is explained graphically in Figure~\ref{fig::delays}.

\begin{lemma}\label{lem::Conv1:AdditiveErrors}
Equation~~\eqref{eq::FP_iter} can be expressed as
\begin{equation}{\label{eq::AsyncInertKMerror}}
x_{k+1} = x_k + \eta (Tx_k-x_k+e_k) = x_k - \eta (Sx_k-e_k)\enspace,
\end{equation}
where $e_k$ is an error term defined as
\begin{equation}{\label{eq::ek}}
e_k := T_A(T_Bx_k + d_k + \beta(c_k-b_k))-T_AT_Bx_k\enspace,
\end{equation}
with
\begin{equation}{\label{eq::a_b_c}}
a_k := \{a_k^i[i]\}_{i=1}^N,\quad b_k := \{b_k^i[i]\}_{i=1}^N, \quad c_k := \{c_k^i[i]\}_{i=1}^N
\end{equation}
and $a_k^i$, $b_k^i$ and $c_k^i$ defined in~\eqref{eq::undo_changes_rw}, and
\begin{equation}{\label{eq::ck}}
 d_k := \gamma Bx_k - \gamma B(x_k-a_k) - a_k \enspace.
\end{equation}
\end{lemma}
The proof of Lemma~\ref{lem::Conv1:AdditiveErrors} is given in Appendix~\ref{app::Conv1:AdditiveError}.

\subsection{Isolate the error. }
The form of the iteration derived in Lemma~\ref{lem::Conv1:AdditiveErrors} will help us separate the error
sequence $\{e_k\}_{k\in\mathbb{N}}$ from the sequence of interest $\{\dist_k\}_{k\in\mathbb{N}}$. To this end, the following result holds, the
derivation of which is given in Appendix~\ref{app::Conv2:IsolateError}.
\begin{lemma}\label{lem::Conv2:IsolateError}
The distance is upper-bounded as:
\begin{equation}{\label{eq::inequality5}}
 \dist_{k+1}^2\leq\left(1-\eta(\nu-\epsilon)\right)\dist_k^2 + \eta\left(\frac{1}{\epsilon}+\frac{\eta(1+\delta)}{\delta}\right)\|e_k\|^2
                                   \enspace,
\end{equation}
for $\eta\in(0, 1/(2(\delta+1)))$ and any $\delta>0$, $\epsilon\in(0,\nu)$, while $\nu$ is the quasi-strong monotonicity constant of $S$ as introduced in Assumption~\ref{Ass::S}.
\end{lemma}
If we manage to bound the last two terms of the equation with
respect to the maximum distance from the set of fixed points, inequality~\eqref{eq::inequality5} will be in the form described by
Lemma~\ref{lem::ISS}. In the next step, we start by bounding the error term $\|e_k\|$.

\subsection{Bound the error recursively. }
We want to bound the error term $\|e_k\|$ by means of $\underset{k-K\leq l\leq k}{\max}\|x_l-x_\ast\|$ for
some $K\in\N$.
We will do so in two phases, first bounding $\|e_k\|$ (given in~\eqref{eq::ek}) recursively with respect to
itself:
\begin{align}{\label{eq::boundek}}
 \|e_k\| \nonumber &= \|T_A(T_Bx_k+d_k+\beta(c_k-b_k))-T_AT_Bx_k\|\\\nonumber
         &\leq \|d_k+\beta(c_k-b_k)\|\\
         &\leq \|d_k\| + \beta\|c_k-b_k\|\enspace,
\end{align}
where the first inequality follows from the nonexpansivity of $T_A$.

It thus suffices to bound $\|d_k\|$ and $\|c_k-b_k\|$ in a recursive way. The result is presented in Lemma~\ref{lem::Conv3:RecursiveBound} below and proven in Appendix~\ref{app::Conv3:RecursiveBound}.

\begin{lemma}\label{lem::Conv3:RecursiveBound}
The quantities $\|a_k\|$, $\|c_k-b_k\|$ and $\|d_k\|$ can be bounded recursively as:
\begin{subequations}
  \begin{alignat}{3}
  \|a_k\| &\leq \eta N\Sigma_{2\tau}(k){\label{eq::boundakRecursive}}\\
  \|c_k-b_k\| &\leq 2\eta N\Sigma_{3\tau}(k){\label{eq::boundbkRecursive}}\\
  \|d_k\| &\leq \eta(1+\gamma L)N\Sigma_{2\tau}(k)\enspace,{\label{eq::boundckRecursive}}
  \end{alignat}
\end{subequations}
where
\begin{equation}{\label{eq::S(K)}}
 \Sigma_K(k):= \sum_{m=k-K}^{k-1}(\|d_m\| + \beta\|c_m-b_m\| + \|Sx_m\|) \tag{$\Sigma$}
 \end{equation}
 and $L$ is the inverse cocoercivity constant of the operator $B$ from Assumption~\ref{Ass::propertiesB}.
\end{lemma}

\subsection{Bound the error with respect to the maximum distance from the set of fixed points. }
Looking at~\eqref{eq::boundakRecursive},~\eqref{eq::boundbkRecursive} and~\eqref{eq::boundckRecursive}, what needs to be bounded is the quantity~\eqref{eq::S(K)},
and consequently the three sums, \ie, $\sum_{m=k-K}^{k-1}\|d_m\|$, $\sum_{m=k-K}^{k-1}\|c_m-b_m\|$ and $\sum_{m=k-K}^{k-1}\|Sx_m\|$ for $K=\{2\tau,3\tau\}$
with respect to the maximum distance from the set of fixed points of $T$.
Lemma~\ref{lem::Conv4:FixedBound} below states the result.

\begin{lemma}\label{lem::Conv4:FixedBound}
 The sequence~\eqref{eq::S(K)} can be upper bounded by the maximum distance from the set of fixed points of $T$ as
 \begin{equation*}
 \|\Sigma_K(k)\| \leq 2K(YN+1)\underset{k-K-3\tau \leq j \leq k-1}{\max}\dist_j\enspace,
\end{equation*}
where $Y:= 1+\gamma L+2\beta$.
Using the above, the error $\|e_k\|$ can be bounded as
\begin{equation}{\label{eq::boundeksquaredExpDist}}
 \|e_k\| \leq \eta X \underset{k-6\tau \leq j \leq k-1}{\max}\dist_j\enspace,
\end{equation}
where
 \[
 X := N(YN+1)(4\tau(1+\gamma L)+6\beta\tau)\enspace.
 \]
\end{lemma}
The Lemma is proven in Appendix~\ref{app::Conv4:FixedBound}.

\subsection{Condition for convergence. }
Let us now recover the condition for the algorithm to converge. By using~\eqref{eq::boundeksquaredExpDist} in~\eqref{eq::inequality5}, we have the desired result expressed as:
\begin{equation*}
 \dist_{k+1}^2 \leq r(\eta)\dist_k^2 + q(\eta)\underset{k-6\tau \leq j \leq k-1}{\max}\dist_j^2 \enspace,
\end{equation*}
where
\begin{equation}{\label{eq::r_and_q}}
 r(\eta):= 1-\eta(\nu-\epsilon), \quad q(\eta):= \eta^3X^2\left(\frac{1}{\epsilon}+\frac{\eta(1+\delta)}{\delta}\right).
\end{equation}

Lemma~\ref{lem::ISS} suggests that \emph{the asynchronous inertial FBS iteration~\eqref{eq::AsyncInertKMerror} will converge to a zero of $S$ at a linear rate $(r(\eta)+q(\eta))^\frac{1}{1+6\tau}$ if the condition}
\begin{equation}{\label{eq::EquationConvergence}}
 1-\eta(\nu-\epsilon) + \eta^3X^2\left(\frac{1}{\epsilon}+\frac{\eta(1+\delta)}{\delta}\right) < 1
\end{equation}
\emph{holds}.
\begin{theorem}{\label{thm::inequality_convergence}}
 Iteration~\eqref{eq::AsyncInertKMerror} will converge at a linear rate as described in~Lemma~\ref{lem::ISS}
 with $r(\eta)$ and $q(\eta)$ given in~\eqref{eq::r_and_q} for
 \[
  \eta < \min\left\{\frac{1}{2(1+\delta)},\frac{1}{X}\sqrt{\frac{2\delta\epsilon(\nu-\epsilon)}{2\delta+\epsilon}}\right\},
 \]
 where $\gamma\in(0,\gamma_{\max})$, $\delta>0$, $\epsilon\in(0,\nu)$ and $\beta>0$. The upper bound $\gamma_{\max}$ ensures that the stepsize $\gamma$ is admissible (a possible option is, \eg, $\gamma_{\max}=2/L$ as proven in Appendix~\ref{app::properties}).
\end{theorem}
Theorem~\ref{thm::inequality_convergence} is proven in Appendix~\ref{app::Conv5:Convergence}.

\section{Connection to other methods}\label{sec::algorithms}
The operators $T_A$ and $T_B$ that constitute the proposed iteration~\eqref{eq::AsyncInertKMerror} give rise to several asynchronous accelerated versions of known algorithms. The form of $T_B$, namely the \emph{forward step} $T_B=I-\gamma B$, played an important role in allowing us to express the errors that arise due to the delays and the inertial term in an additive manner. At the same time, this limits the applicability of our iteration to methods that can be cast as forward-backward iterations. Below, we introduce the extensions of some popular algorithms that can be seen as special cases of iteration~\eqref{eq::AsyncInertKMerror}.

\subsection{Gradient descent. }{\label{subsec::gradient}}
 Classical gradient descent can be recovered by choosing $T_A=I$ and $T_B=I-\gamma\nabla f$, for a differentiable strongly convex function $f:\reals^n\mapsto\reals$ with Lipschitz continuous gradient. The inertial asynchronous iteration becomes:
 \[
   x_{k+1}[i] = (1-\eta)x_{k}[i] + \eta\underbrace{\left(x_\mathrm{read}^i[i]-\gamma\nabla_{i}f(x_\mathrm{read}^i) + \beta(y_\mathrm{write}-y_\mathrm{write}^\mathrm{prev})\right)}_{z^i}\enspace,
 \]
 where $\nabla_i f(x) := \nabla_{x_i} f(x): \reals^n\mapsto\reals^{n_i}$, $\sum_{i=1}^Nn_i=n$, and corresponds to an asynchronous iteration of the Heavy Ball gradient method.

 \subsection{Proximal gradient. }{\label{subsec::prox_gradient}}
 Let us consider the optimization problem
 \begin{equation}{\label{eq::FBS}}
  \mbox{minimize}\; f(x) + \sum_{i=1}^N g_i(x_i)\enspace,
 \end{equation}
 where $f$ is stronlgy convex differentiable with Lipschitz continuous gradient, while $g_i\in\Gamma_0(\reals^{n_i})$ and $x\in\reals^n=\reals^{n_1}\times\cdots\times\reals^{n_N}$.
 The iteration reads:
  \[
   x_{k+1}[i] = (1-\eta)x_{k}[i] + \eta\underbrace{\prox_{\gamma g_{i}}\left(x_\mathrm{read}^i[i]-\gamma\nabla_{i}f(x_\mathrm{read}^i)+\beta\left(y_\mathrm{write}-y_\mathrm{write}^\mathrm{prev}\right)\right)}_{z^i}\enspace,
 \]
and corresponds to a relaxed and asynchronous version of the proximal Heavy Ball method in~\cite{ipiasco}.

\subsection{Other methods. }{\label{subsec::others}}
Besides the two instances analyzed above, a variety of convex optimization algorithms can be expressed as a forward-backward iteration, and consequently give rise to novel asynchronous implementations. The \emph{generalized forward backward splitting}~\cite{GFBS,PrecGFBS}, the \emph{forward-Douglas-Rachford splitting}~\cite{FDRS,FDRS_rate} and several \emph{primal-dual optimization methods}~\cite{combettes2014forward} can be viewed as candidates, just to name a few.

\section{Application: Distribution network real-time dispatch }\label{sec::example}
We consider the problem of cooperative tracking of a reference signal from a population of controllable buildings in combination with a battery energy storage system. These problems typically arise in the context of microgrids, where a mixture of energy generation, energy storage elements and loads are coupled together in order to satisfy a predicted power demand profile.

The growing interest in turning once passive loads into active prosumers, along with the introduction of dispatchable energy storage elements in the grid, is motivated by a rapid and significant increase of renewable production into the generation mix. Renewable energy sources are inherently uncertain and volatile, and they pose new challenges to the classic control paradigm of the power grid. At the same time, the introduction of demand side management via large loads as, \eg, commercial buildings, raises privacy concerns and calls for distributed implementations. Several paradigms have emerged toward this direction, see, \eg,~\cite{EPFL-ARTICLE-222513}.

Our goal is to track a 15-minute resolution trajectory, called the dispatch plan that is computed one day before the beginning of operation. This is achieved by modulating the power consumption of a grid-connected battery energy storage system (BESS) and of the thermal consumption of a fleet of commercial controllable buildings (CB). The problem has been proposed and solved in~\cite{EPFL-CONF-226792} using as benchmark an experimental setup with one controllable office and a large battery. In this example we scale up the problem by considering several CB's and we compare synchronous versus asynchronous implementations, including our proposed scheme.

\subsection{Modeling the agents. }
The grid comprises the following entities:
\paragraph{Controllable Loads:} Small, medium and large office buildings, generated by the OpenBuild software~\cite{OpenBuild}. The buildings are described as linear dynamical systems,
the input to which is the thermal heat $(kW)$ that is entering or leaving each zone, while the output is the temperature at each zone $({}^\circ C)$.
The energy conversion systems (electrical to thermal) is modeled as a static map, which is represented by a constant coefficient of performance (COP).
The buildings can be dispatchable by \emph{increasing or decreasing their consumption with respect to some baseline power profile}.
An individual building operates within temperature constraints. The local optimization problem for the CB $i$ becomes
\begin{equation}{\label{eq::BuildingProgram}}
g^\mathrm{cb}_i(p^{\mathrm{cb}}_i,u_i,x_i,y_i) := \Big\{\frac{1}{2}\|y_i(t)-T_i^\mathrm{ref}(t)\|^2_2
                                             \;\mid\;(p^{\mathrm{cb}}_i,u_i,x_i,y_i)\in\mathcal{C}^\mathrm{cb}_i \Big\}\enspace,
\end{equation}
\[
\mathcal{C}^\mathrm{cb}_i =
\left \{
                  \begin{array}{ll}
                    x_i(t+1) = A_ix_i(t)+B_{u,i}u_i(t)+B_{w,i}\hat{w}_i(t) \\
                    x_i(0) = x_i^\mathrm{init} \\
                    y_i(t) =C_ix_i(t) \\
                    y_{\min,i}(t)\leq y_i(t)\leq y_{\max,i}(t) \\
                    u_{\min,i}\leq u_i(t)\leq u_{\max,i} \\
                    p^{\mathrm{cb}}_i(t) = \sum_{j=1}^{N_i}u_{ij}(t) \\
                  \end{array}
\right\}\enspace,
\]
with $x_i\in\reals^{n_iT},u_i=\{u_{ij}\}_{j=1}^{m_i}\in\reals^{m_iT}$, where $p^{\mathrm{cb}}_i(t)$ is the total amount (electrical equivalent)
of the thermal consumption at time $t$, and $m_i$ is the number of zones of the building. Zone temperatures are described with the variables $y_i\in\reals^{l_iT}$.
In addition, the temperature constraints are relaxed outside working hours, hence the time varying constraint limits $y_{\min,i}(t), y_{\max,i}(t)$
(see Table~\ref{table::tab_micro-grid}). The desired zone temperature is denoted with $T_i^\mathrm{ref}$.
\paragraph{Storage:} The setup is completed with a grid-connected Lithium Titanate grid-connected $500\mathrm{kWh}$ BESS.
The battery is represented as single-state linear model, with the state being the state-of-charge (SOC) and the input being the active power denoted by $p^\mathrm{bess}(t)\in\reals$. The battery operates within capacity and power limits, while the purpose is to keep the SOC close to a reference value $SOC^\mathrm{ref}(t)\in\reals$.
\begin{equation}{\label{eq::BESSProgram}}
g^\mathrm{bess}(p^\mathrm{bess}) := \Big\{\frac{1}{2}\sum_{t=1}^T\|SOC(t)-SOC^\mathrm{ref}(t)\|^2_2
                                             \;\mid\;p^\mathrm{bess}\in\mathcal{C}^\mathrm{bess} \Big\}\enspace,
\end{equation}
\[
\mathcal{C}^\mathrm{bess} =
\left \{
                  \begin{array}{ll}
                    SOC(t+1) = a SOC(t) + b p^\mathrm{bess}(t) \\
                    SOC(0) = SOC^\mathrm{init} \\
                    SOC_{\min}  \leq SOC(t) \leq SOC_{\max} \\
                    p_{\min}^\mathrm{bess} \leq p^\mathrm{bess}(t) \leq p_{\max}^\mathrm{bess}\\
                  \end{array}
\right\}\enspace.
\]

\begin{table}[H]
\centering
\scalebox{0.5}{
\Large
\begin{tabular}{| l l l l c |}
\hline
& & & &\\
\multicolumn{5}{| l |}{\textbf{Simulation characteristics}} \\
& & & &\\
Data & \multicolumn{3}{ l }{ $1^{\mathrm{st}}$ January 2000} & \\
Location & \multicolumn{3}{ l }{ Lausanne} & \\
Time & \multicolumn{3}{ l }{00:00 - 24:00} & \\
Sampling time & \multicolumn{3}{ l }{15} & $\mathrm{min}$ \\
Horizon & \multicolumn{3}{ l }{96} & $\mathrm{-}$ \\
& & & &\\

\multicolumn{5}{| l |}{\textbf{Buildings}} \\
& & & &\\
Minimum temperature (day/night) & \multicolumn{3}{ l }{20/18} & $\mathrm{^\circ C}$ \\
Maximum temperature (day/night) & \multicolumn{3}{ l }{24/28} & $\mathrm{^\circ C}$ \\
Heat pump $\mathrm{COP}$ & \multicolumn{3}{ l }{3.0} & $\mathrm{-}$ \\
& & & &\\

\multicolumn{5}{| l |}{} \\
& \textbf{Small}&\textbf{Medium}&\textbf{Large} &\\
Number of systems (Case A, B, C, D) & 3/6/14/32 & 2/4/5/16 & 0/0/1/2 & $\mathrm{\#}$ \\
Area & 511 & 4982 & 46320 & $\mathrm{m^2}$ \\
Tariff (day/night) & 21.6/12.7 & 13.15/8.3 & 13.15/8.3 & $\mathrm{ct./kWh}$ \\
Number of states & 15 & 54 & 57 & $\mathrm{-}$ \\
Number of inputs & 5 & 18 & 19 & $\mathrm{-}$ \\
Average thermal consumption & 4 & 40 & 75 & $\mathrm{W/m^2}$ \\
Average computation time (prox per agent) & $0.070\pm 0.010$ & $0.243\pm 0.005$ & $0.267\pm 0.001$ & $\mathrm{sec}$ \\
& & & &\\

\multicolumn{5}{| l |}{\textbf{Battery}} \\
Energy storage capacity & \multicolumn{3}{ l }{500} & $\mathrm{kWh}$ \\
C-rate & \multicolumn{3}{ l }{0.2} & $\mathrm{-}$ \\
Material & \multicolumn{3}{ l }{Lithium-ion} &\\
Average computation time (prox) & \multicolumn{3}{ l }{$0.023\pm 0.003$} & $\mathrm{sec}$ \\
& & & &\\

\hline
\end{tabular}}
\caption{Micro-grid case study overview}
\label{table::tab_micro-grid}
\end{table}

\subsection{Modeling the dispatch problem. }
The dispatch problem can be cast as:
\begin{subequations}
\begin{align}
&{\text{minimize}}   && \frac{1}{2}\sum_{t=1}^T\Big(\textcolor{NavyBlue}{\|SOC(t)-SOC^\mathrm{ref}(t)\|^2_2} + \textcolor{Red}{\alpha_1\|p^\mathrm{bess}(t)\|^2_2}\Big) {\label{eq::prob_dispatch::SOC}} \\
&                    && + \frac{1}{2}\sum_{i=1}^N\sum_{t=1}^T \Big(\textcolor{ForestGreen}{\|y_i(t)-T_i^\mathrm{ref}(t)\|^2_2} + \textcolor{Red}{\alpha_1\|p_i^{\mathrm{cb}}(t)-\hat{p}_i^\mathrm{cb}(t)\|_2^2} \Big) {\label{eq::prob_dispatch::regular}} \\
&                    && + \textcolor{Red}{\frac{\alpha_2}{2}\sum_{t=1}^T\|p^\mathrm{bess}(t) + \sum_{i=1}^N(p_i^\mathrm{cb}(t)-\hat{p}_i^\mathrm{cb}(t))-r(t)\|^2_2} {\label{eq::prob_dispatch::consensus}}\\
&{\text{subject to}} && \textcolor{ForestGreen}{(p^{\mathrm{cb}}_i,u_i,x_i,y_i)\in\mathcal{C}^\mathrm{cb}_i}, \; i=1,\ldots,N \\
&                    && \textcolor{NavyBlue}{p^\mathrm{bess}\in\mathcal{C}^\mathrm{bess}}\enspace,
\end{align}
\label{eq::prob_dispatch}
\end{subequations}
with variables $p^{\mathrm{cb}}_i,\;i=1,\ldots,N$ and $p^\mathrm{bess}$, and the variables $u_i,x_i,y_i$ local to CB $i$.

Equations~\eqref{eq::prob_dispatch::SOC}~and~\eqref{eq::prob_dispatch::regular} express the deviation of the BESS SOC from its reference value, set to $SOC^\mathrm{ref}(t)=0.8SOC_{\max}$, as well as the deviation of the indoor temperature from its reference value (see~Table~\ref{table::tab_micro-grid}). The additional quadratic terms penalized with $\alpha_1=10^{-2}$ are introduced for regularization purposes. Equation~\eqref{eq::prob_dispatch::consensus} expresses \emph{the deviation of the aggregate buildings' flexibility $\sum_{i=1}^N(p_i^\mathrm{cb}(t)-\hat{p}_i^\mathrm{cb}(t))$ along with the BESS flexibility $p^\mathrm{bess}$ from the given reference $r$.} In a perfectly dispatchable network, this term should be put to zero, hence it is penalized much more heavily than the other terms with $\alpha_2=10^4$.

\subsection{Simulation setup. }
Our purpose is to solve~\eqref{eq::prob_dispatch} by means of the synchronous, the asynchronous and the inertial asynchronous versions of the FBS algorithm.
To this end, let us make the problem more compact by grouping the terms. The terms depicted in blue color are private to the BESS system, the terms in green are private to the CB agents, while red terms comprise the global objective, denoted hereafter as $f(p^\mathrm{bess},p^\mathrm{cb})$. Note that the local subproblems in blue and green correspond to the quadratic programs (QP)~\eqref{eq::BESSProgram} and~\eqref{eq::BuildingProgram}, respectively. Since each variable $p^\mathrm{cb}_i$ is private to agent $i$ and $f$ couples all the variables through the (strongly convex) quadratic objective, the problem takes the form~\eqref{eq::FBS} and is consequently solved using the proximal gradient method.

We consider four case studies (A, B, C and D), namely a mix of the BESS and $N=5,10,20$ and $50$ CB's. The tracking signal $r$ that we assume is the realized
\emph{Area Control Signal (ACS)}, as it was broadcast by the Swiss grid operator~\cite{ACS}, for
the $1^{\mathrm{st}}$ of January of the year 2000, scaled down by the appropriate factor in each of the four cases so that it becomes (almost) trackable by our mix. The prediction horizon has a length of $24$ hours, or $T=96$ in $15$ minutes intervals.

The source of asynchronicity in this framework is the diverse computational load of the different agents (CB's and the BESS). Although problems \eqref{eq::BuildingProgram} and \eqref{eq::BESSProgram} are QP's, their size varies greatly with the number of states and inputs, as indicated in Table~\ref{table::tab_micro-grid}. The delay $\tau$ is, therefore, computed based on the number of the updates per agent in a unit of time. In order to compute this number, we solve 100 proximal minimization steps per agent and fit a normal distribution to the solve times. The average computation time is then used to decide upon the frequency of the updates. The communication delays are assumed to be zero in the simulation. The proximal minimization problems are solved using the YALMIP optimizer~\cite{YALMIP} with the Gurobi solver. Finally, the relaxation parameter is set to $\eta=0.9$, which, in spite of the (much) smaller value suggested by Theorem~\ref{thm::inequality_convergence}, worked well in our setting.

Problem~\eqref{eq::prob_dispatch} is solved using the proximal gradient method. A comparison between (i) the synchronous version of the method (all CB's and BESS update before a new gradient $\nabla f$ is communicated), (ii) the asynchronous version with coordinate updates (only $i=i_k$ updates at each gobal clock count, with $i=1,\ldots,N+1$ and $N+1$ corresponds to the BESS agent), (iii) the asynchronous aggregated version (\eqref{eq::FBS} with $\beta=0$) and (iv) the asynchronous inertial aggregated version (\eqref{eq::FBS} with $\beta=0.99$). Table~\ref{table::accuracy} depicts the accuracy reached within $Ts=40\mathrm{sec}$ of simulated wall-clock time using the four algorithms presented above in the four case studies. Table~\ref{table::updates} presents the average number of updates per type of agent within these $40\mathrm{sec}$.

\begin{figure}[!htb]
  \begin{center}
  \includegraphics[scale=0.5]{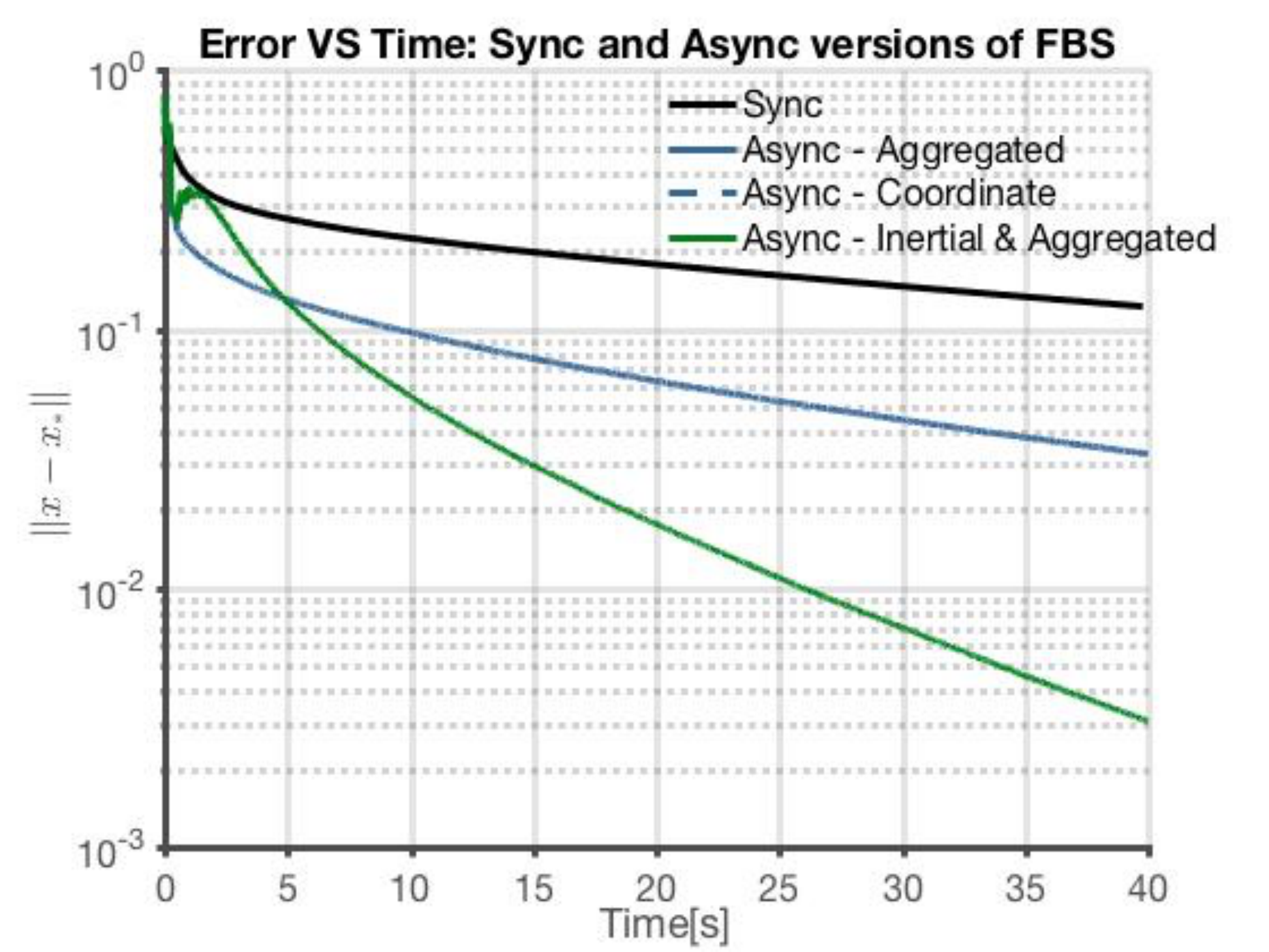}
  \caption{Distance from optimizer VS wall-clock time.}
  \label{fig::convergence}
  \end{center}
\end{figure}

\begin{figure}[!htb]
  \begin{center}
  \includegraphics[scale=0.5]{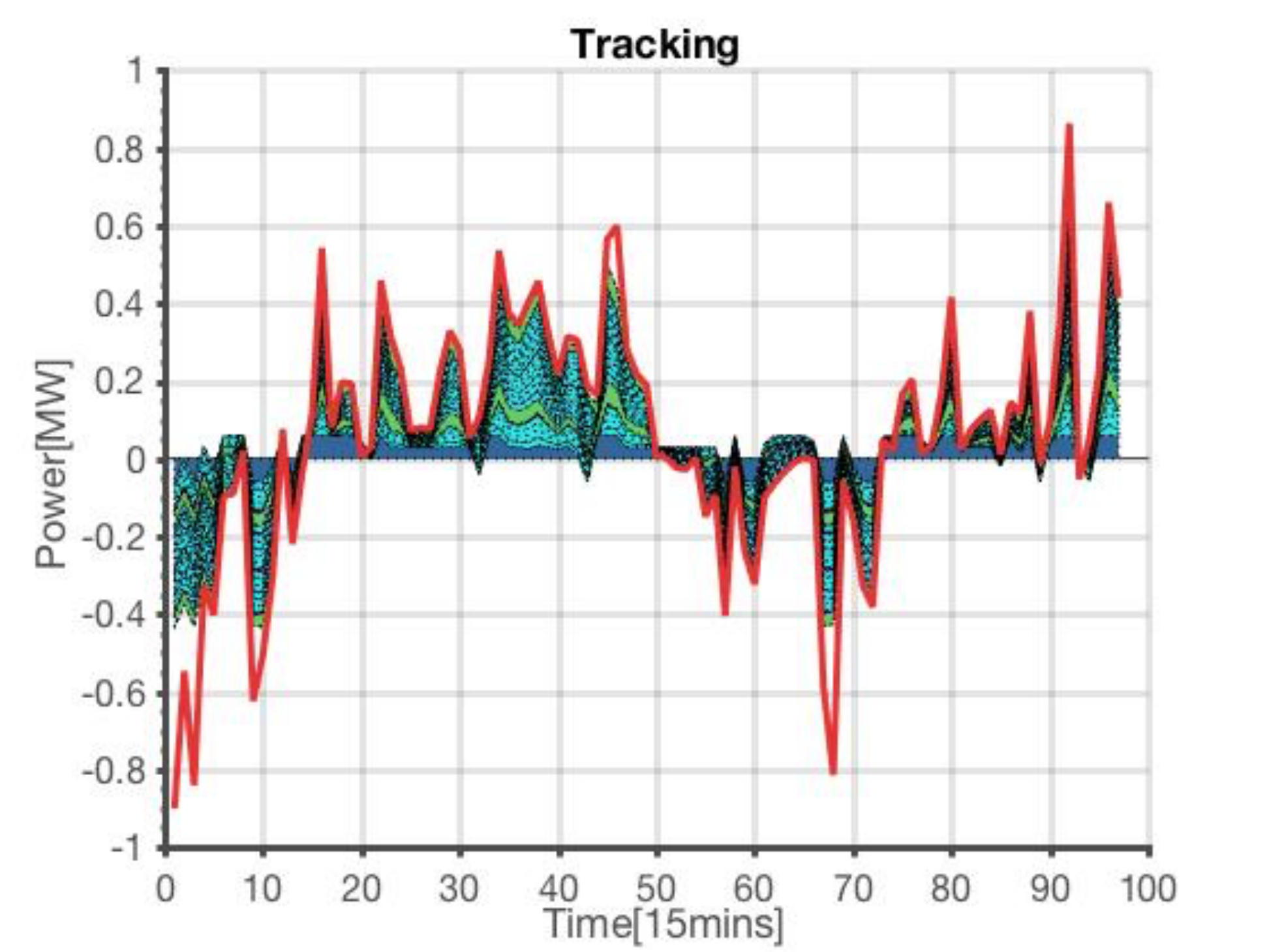}
  \caption{The reference signal to be tracked is depicted in red. The contribution of the BESS is colored in dark blue, that of the medium scale buildings in lighter blue and that of the two large buildings in green. The contribution of the small buildings is colored in pink, but is hardly visible due to their small capacity and despite their large population.}
  \label{fig::area}
  \end{center}
\end{figure}

Several conclusions can be derived. First, the asynchronous version with coordinate updates and the asynchronous aggregated version of the proximal gradient method are almost identical in performace, thus there is neither deterioration (at least in the simulated cases) nor improvement when using the old updates. Second, the asynchronous versions perform considerably better than their synchronous counterpart in all cases. This is an expected outcome since the larger the load imbalance among the agents, the more the algorithm benefits from the asynchronicity, as suggested by the number of updates per agent in Table~\ref{table::updates}. Finally, the proposed inertial acceleration scheme results in considerably better performance in terms of speed of convergence in all cases. A graphical depiction of the convergence performance of the four methods for $N=5$ is given in Figure~\ref{fig::convergence}. For Case D ($N=50$), the area plot in Figure~\ref{fig::area} elaborates on the contribution of each of the agents in tracking the reference signal.

\begin{table}
  \normalsize
  \begin{center}
  \scalebox{1.0}{
  \begin{tabular}{|l||*{4}{c|}}
     \cline{1-5}
     \backslashbox{$\textbf{N}$}{\textbf{Algo.}} & \thead{\textbf{Sync}} & \thead{\textbf{Async} \\ \textbf{Coordinate}} & \thead{\textbf{Async} \\ \textbf{Aggregated}} & \thead{\textbf{Async} \\ \textbf{Agg. Inert.}} \\ \cline{1-5}
       5  & 0.116 & 0.030 & 0.030 & 0.003 \\ \cline{1-5}
       10 & 0.252 & 0.061 & 0.061 & 0.012 \\ \cline{1-5}
       20 & 0.315 & 0.078 & 0.078 & 0.015 \\ \cline{1-5}
       50 & 0.824 & 0.649 & 0.649 & 0.448\\ \cline{1-5}
  \end{tabular}
   }
\end{center}
\caption{Accuracy reached within $Ts=40\mathrm{sec}$.}
{\label{table::accuracy}}
\end{table}

\begin{table}
  \normalsize
  \begin{center}
  \scalebox{0.7}{
  \begin{tabular}{|l||*{4}{c|}}
     \cline{1-5}
     \backslashbox{\textbf{Algo.}}{$\textbf{N}$} & 5            & 10   &  20     & 50 \\ \cline{1-5}
       \textbf{Sync}                             & 172          & 156  &  148    & 137 \\ \cline{1-5}
       \textbf{Async} \par \textbf{Agg. Inert.}  & $\begin{matrix} 1639 \\ 565.6/171/-\end{matrix}$  & $\begin{matrix} 1798 \\ 588.5/157/-\end{matrix}$ & $\begin{matrix} 1674 \\ 515.6/163.2/147\end{matrix}$ & $\begin{matrix} 1101 \\ 426.37/148.50/142.50\end{matrix}$\\ \cline{1-5}
  \end{tabular}
   }
\end{center}
\caption{Average number of updates per agent within simulation time.}
{\label{table::updates}}
\end{table}

\section{Conclusions}\label{sec::conclusions}
We proposed an inertial and asynchronous forward-backward iteration for solving monotone inclusion problems. The iteration is tailored for distributed convex optimization problems and differs from existing approaches since (i) the component updates are selected in a deterministic way, (ii) older updates contribute to the upcoming one in a fashion resembling aggregated gradient methods and (iii) the iteration hosts a momentum term that speeds up the practical convergence. We derived new versions of three commonly used methods stemming from our approach, and we illustrated the effectiveness of the method when used to solve an optimal dispatch problem in a distribution grid with a pool of heterogeneous energy resources.

There is plenty of space for improving the proposed approach and the like. The first things that naturally come to mind regard dropping the strong monotonicity (strong convexity) assumption and using an optimal Nesterov-like momentum sequence instead of a fixed scalar value. In addition, although the momentum sequence practically boosts the performance of the scheme, the current convergence analysis does not exhibit its benefits. On the contrary, our analysis treats the additional degree of freedom that momentum offers as an extra perturbation. Since it is known that inertial acceleration improves the rate at which the sequence of iterates converges to a solution in the strongly convex case, we suspect that a similar result is applicable to the asynchronous framework.

\section*{Acknowledgments}
The research leading to these results has received funding from the European Research Council under the European Union's Seventh Framework Programme (FP/2007-2013) / ERC Grant Agreement n. 307608: BuildNet.

\appendix  

\section{Proof of Lemma~\ref{lem::Conv1:AdditiveErrors} }{\label{app::Conv1:AdditiveError}}
{\it Proof}  Using Assumption~\ref{Ass::propertiesB} and equations~\eqref{eq::undo_changes_rw}, we can start by rewriting $y_B=(T_B x^i_\mathrm{read})[i]$.
\begin{align*}
 (T_B x^i_\mathrm{read})[i] &= x^i_\mathrm{read}[i] - \gamma (Bx^i_\mathrm{read})[i] \\
                     &= x_{k}[i] - a_k^i[i] - \gamma (Bx_{k})[i] + \gamma (Bx_{k})[i] - \gamma (Bx^i_\mathrm{read})[i] \\
                     &= (T_Bx_k)[i] + \gamma \left((Bx_{k})[i]-(Bx^i_\mathrm{read})[i]\right) - a_k^i[i] \\
                     &= (T_Bx_k)[i] + d_k^i[i]\enspace,
\end{align*}
where $d_k^i[i]=\gamma \left((Bx_{k})[i]-(Bx^i_\mathrm{read})[i]\right) - a_k^i[i]$.

Similarly, we have from equations~\eqref{eq::undo_changes_rw} that
\begin{align*}
 \beta ( y_\mathrm{write} - y_\mathrm{write}^\mathrm{prev} ) &= \beta ( x_k[i] - b_k^i[i] - x_k[i] + c_k^i[i] ) \\
                     &= \beta (c_k^i[i] - b_k^i[i])\enspace.
\end{align*}

Using the above relations, a coordinate update of iteration~\eqref{eq::FP_iter} can be expressed as
\begin{equation*}
x_{k+1}[i] = x_k[i] + \eta \Big(T_{A_{i}} \left( (T_B x_k)[i] + d_k^i[i] + \beta (c_k^i[i] - b_k^i[i]) \right) - x_k[i] \Big)\enspace,
\end{equation*}
or, equivalently, as
\begin{equation*}
x_{k+1}[i] = x_k[i] + \eta \left( T_{A_{i}}((T_B x_k)[i]) - x_k[i] + e_k[i]\right)\enspace,
\end{equation*}
with $e_k[i] = T_{A_{i}}\left((T_Bx_k)[i] + d_k^i[i] + \beta(c_k^i[i]-b_k^i[i])\right) - T_{A_{i}}((T_Bx_k)[i])$, which concludes the proof.
\qed

\section{Proof of Lemma~\ref{lem::Conv2:IsolateError} }{\label{app::Conv2:IsolateError}}
{\it Proof} Squaring~\eqref{eq::AsyncInertKMerror} we get:
\begin{equation}{\label{eq::inequality}}
 \|x_{k+1}-x_\ast\|^2 = \|x_{k}-x_\ast\|^2 -2\eta\langle x_k-x_\ast, Sx_k-e_k \rangle +\eta^2\|Sx_k-e_k\|^2\enspace.
\end{equation}

Let us now analyze the second and third term in~\eqref{eq::inequality}.
\begin{itemize}
 \item Bound $-2\eta\langle x_k-x_\ast, Sx_k-e_k \rangle$: We will upper-bound the resulting inner product terms.
 In order to do so, we use both the cocoercivity and the quasi-strong monotonicity of $S$,
 the former proven in Appendix~\ref{app::properties}, and the latter holding from Assumption~\ref{Ass::propertiesB}.
Since $S$ is $1/2$-cocoercive, we have that
\[
 \langle x_k-x_\ast,Sx_k\rangle \geq \frac{1}{2}\|Sx_k\|^2.
\]
From the quasi-$\nu$-strong monotonicity of $S$ we have:
\[
 \langle x_k-x_\ast,Sx_k\rangle \geq \nu\|x_k-x_\ast\|^2.
\]
Putting these two together, we get that
\begin{equation}{\label{eq::inequality3}}
 -2\eta\langle x_k-x_\ast,Sx_k\rangle \leq -\eta\nu\dist_k^2
                   -\frac{\eta}{2}\|Sx_k\|^2\enspace.
\end{equation}
For the second inner product term involving the error we can easily derive the bound
\begin{equation}{\label{eq::inequality_product1}}
 2\eta\langle x_k-x_\ast,e_k\rangle \leq 2\eta\dist_k\|e_k\|\enspace.
\end{equation}
Equations~\eqref{eq::inequality3} and~\eqref{eq::inequality_product1} result in the bound
\begin{equation}{\label{eq::inner_product1}}
-2\eta\langle x_k-x_\ast, Sx_k-e_k\rangle \leq -\eta\nu\dist_k^2-\frac{\eta}{2}\|Sx_k\|^2 + 2\eta\dist_k\|e_k\|\enspace.
\end{equation}

\item Bound $\eta^2\|Sx_k-e_k\|^2$: By developing the square, we have that
\begin{align}{\label{eq::inequality_product2}}
\eta^2\|Sx_k-e_k\|^2 &= \eta^2 \left(\|Sx_k\|^2 -2\langle Sx_k,e_k\rangle + \|e_k\|^2\right)\enspace.
\end{align}
The inner product term in~\eqref{eq::inequality_product2} can be bounded by employing Young's inequality\footnote{For two nonnegative real numbers $x$ and $y$, it holds that $xy\leq \frac{\delta x^2}{2}+\frac{y^2}{2\delta}$ for every $\delta>0$.} as follows:
\begin{align}{\label{eq::Ek_Young}}
-2\langle Sx_k,e_k\rangle \nonumber&\leq 2\|Sx_k\|\|e_k\| \\\nonumber
     &\leq 2(\frac{\delta}{2}\|Sx_k\|^2+\frac{1}{2\delta}\|e_k\|^2)\\
     &= \delta\|Sx_k\|^2 + \frac{1}{\delta}\|e_k\|^2\enspace,
\end{align}
for any $\delta>0$. Putting together~\eqref{eq::inequality_product2} and~\eqref{eq::Ek_Young}, we get the bound:
\begin{equation}{\label{eq::inner_product2}}
\eta^2\|Sx_k-e_k\|^2 \leq \eta^2(1+\delta)\|Sx_k\|^2 +
\eta^2\frac{(\delta+1)}{\delta}\|e_k\|^2\enspace.
\end{equation}
\end{itemize}

Using~\eqref{eq::inner_product1} and~\eqref{eq::inner_product2}, inequality~\eqref{eq::inequality} can be written as
\begin{equation}{\label{eq::inequality4}}
 \dist_{k+1}^2 \leq (1-\eta\nu)\dist_k^2
                            +\eta\left(-\frac{1}{2}+\eta(1+\delta)\right)\|Sx_k\|^2
                            +2\eta\dist_k\|e_k\|
                            +\eta^2\frac{(\delta+1)}{\delta}\|e_k\|^2
                                   \enspace.
\end{equation}

The second term in the sum can be eliminated by asumming that
\begin{equation}{\label{eq::rhoFirstBound}}
 -\frac{1}{2}+\eta(1+\delta)<0\Rightarrow \eta<\frac{1}{2(1+\delta)}\enspace,
\end{equation}
which gives rise to the inequality
\begin{equation}{\label{eq::inequality6}}
 \dist_{k+1}^2 \leq (1-\eta\nu)\dist_k^2
                            +2\eta\dist_k\|e_k\|
                            +\eta^2\frac{(\delta+1)}{\delta}\|e_k\|^2
                                   \enspace.
\end{equation}
The complicating term on the right hand side can be eliminated by using once more Young's inequality, \ie,
\begin{align*}
 2\eta\dist_k\|e_k\| &
                     \leq 2\eta\left(\frac{\epsilon}{2}\dist_k^2+\frac{1}{2\epsilon}\|e_k\|^2\right)\\
                               &= \eta\epsilon\dist_k^2+\frac{\eta}{\epsilon}\|e_k\|^2\enspace.
\end{align*}
\qed

\section{Proof of Lemma~\ref{lem::Conv3:RecursiveBound} }{\label{app::Conv3:RecursiveBound}}
{\it Proof}
\begin{itemize}
  \item Let us start with bounding $\|a_k\|$ for some arbitrary $k\in\N_{+}$, from the definition of which (eq.\eqref{eq::a_b_c}) we have:
\begin{align}{\label{eq::boundakRecursive1}}
 \|a_k\| \nonumber &\leq N\underset{1\leq i\leq N}{\max}\|a_k^i[i]\|\\\nonumber
                   &= N\|\sum_{m=k-l_{i_\ast}}^{k-1}(x_{m+1}[i_\ast]-x_m[i_\ast])\|\\\nonumber
                   &\leq N\sum_{m=k-l_{i_\ast}}^{k-1}\|x_{m+1}[i_\ast]-x_m[i_\ast]\|\\\nonumber
                   &\leq N\sum_{m=k-l_{i_\ast}}^{k-1}\|x_{m+1}-x_m\|\\\nonumber
                   &\leq N\sum_{m=k-2\tau}^{k-1}\|x_{m+1}-x_m\|\\\nonumber
                   &= \eta N\sum_{m=k-2\tau}^{k-1}\|e_m-Sx_m\|\\\nonumber
                   &\leq \eta N\left(\sum_{m=k-2\tau}^{k-1}\|e_m\|+\sum_{m=k-2\tau}^{k-1}\|Sx_m\|\right)\\
                   &\leq \eta N\sum_{m=k-2\tau}^{k-1}(\|d_m\| + \beta\|c_m-b_m\| + \|Sx_m\|)\enspace,
\end{align}
where we denoted as $i_\ast=\underset{i\in\{1,\ldots,N\}}{\argmax}\|x_k[i]-x_{k-l_i}[i]\|$. The fourth inequality follows from the definitions of $a_k$ in~\eqref{eq::undo_changes_rw}, the second equality from~\eqref{eq::AsyncInertKMerror}, while the last two inequalities from the triangle inequality and~\eqref{eq::boundek}.

 \item Upper-bounding $\|d_k\|$ can be achieved by using the $1/L$-cocoercivity of $B$ (Assumption~\ref{Ass::propertiesB}) in~\eqref{eq::ck}:
\begin{equation}{\label{eq::boundcdtemp}}
 \|d_k\|\leq (1+\gamma L)\|a_k\|\enspace.
\end{equation}
To obtain a recursive bound for $\|d_k\|$, we need a recursive bound for $\|a_k\|$, which has been already computed from~\eqref{eq::boundakRecursive1}. We thus just need to derive a similar one for the difference $\|c_k-b_k\|$.

 \item Bound $\|c_k-b_k\|$: Following the exact same process as in~\eqref{eq::boundakRecursive1}, we have that
  \begin{align}{\label{eq::bound_ad_diff}}
  \|c_k-b_k\| \nonumber&\leq \|c_k\| + \|b_k\| \\
              \nonumber&\leq N\sum_{m=k-3\tau}^{k-1}\|x_{m+1}-x_m\| + N\sum_{m=k-2\tau}^{k-1}\|x_{m+1}-x_m\|\\\nonumber
                                     &\leq \eta N\Big(\sum_{m=k-3\tau}^{k-1}(\|d_m\| + \beta\|c_m-b_m\| + \|Sx_m\|) +\\
                                     &\quad\sum_{m=k-2\tau}^{k-1}(\|d_m\| + \beta\|c_m-b_m\| + \|Sx_m\|)\Big)\enspace.
 \end{align}
\end{itemize}
Finally,~\eqref{eq::boundakRecursive1},~and~\eqref{eq::bound_ad_diff} can be bounded by means of the quantity~\eqref{eq::S(K)}, and by substituting \eqref{eq::boundakRecursive1} to~\eqref{eq::boundcdtemp}, the result follows. \qed

\section{Proof of Lemma~\ref{lem::Conv4:FixedBound} }{\label{app::Conv4:FixedBound}}
{\it Proof}
Let us start by bounding the quantities involved in~\eqref{eq::S(K)}, namely $\|a_k\|$, $\|b_k\|$, $\|c_k\|$ and $\|d_k\|$
with respect to the maximum distance from the optimizer. Note that $\forall\;i=1,\ldots,N$ and for $l_i\in\{1,\ldots,3\tau\}$ holds
\begin{align*}
\|x_k[i]-x_{k-l_i}[i]\|  & \leq \|x_k-x_{k-l_i}\| \\
                         & = \|x_k-x_\ast+x_\ast-x_{k-l_i}\| \\
                         & \leq \|x_k-x_\ast\| + \|x_\ast-x_{k-l_i}\|\enspace.
\end{align*}
Since the first inequality holds $\forall\;i=1,\ldots,N$, by denoting $i_\ast=\underset{i\in\{1,\ldots,N\}}{\argmax}\|x_k[i]-x_{k-l_i}[i]\|$ and subsequently taking the $\max$ among all $i$'s and multiplying both sides by $N$, we get
\begin{align*}
N\underset{1\leq i\leq N}{\max}\|x_k[i]-x_{k-l_i}[i]\| &= N\|x_k[i_\ast]-x_{k-l_{i_\ast}}[i_\ast]\| \\
                                                       &\leq N\left(\|x_k-x_\ast\| + \|x_\ast-x_{k-l}\|\right) \enspace,
\end{align*}
where $l=l_{i_\ast}$.

From the definition of $a_k$ in~\eqref{eq::undo_changes_rw} we have that
\begin{align}{\label{eq::boundakDist1}}
\|a_k\| \nonumber & \leq N\underset{1\leq i\leq N}{\max}\|x_k[i]-x_{k-l_i}[i]\| \\
        \nonumber & \leq N\left(\|x_k-x_\ast\| + \|x_\ast-x_{k-l}\|\right) \\
        \nonumber & \leq N\left(\|x_k-x_\ast\| + \underset{k-2\tau \leq m \leq k-1}{\max}\dist_m\right) \\
                  & \leq 2N\underset{k-2\tau \leq m \leq k}{\max}\dist_m\enspace.
\end{align}

Following developments similar to~\eqref{eq::boundakDist1}, we conclude that
\begin{align}{\label{eq::boundDist}}
\|a_k\| \nonumber & \leq 2N\underset{k-2\tau \leq m \leq k}{\max}\dist_m \\
\|b_k\| \nonumber & \leq 2N\underset{k-2\tau \leq m \leq k}{\max}\dist_m \\
\|c_k\| \nonumber & \leq 2N\underset{k-3\tau \leq m \leq k}{\max}\dist_m \\
\|d_k\|           & \leq 2(1+\gamma L)N\underset{k-2\tau \leq m \leq k}{\max}\dist_m
\end{align}

Using~\eqref{eq::boundDist}, the sums can be easily bounded as shown below.
\begin{align}{\label{eq::boundSumDist}}
\sum_{m=k-K}^{k-1}\|a_m\| \nonumber & \leq 2NK\underset{k-K-2\tau \leq j \leq k-1}{\max}\dist_j \\
\sum_{m=k-K}^{k-1}\|b_m\| \nonumber & \leq 2NK\underset{k-K-2\tau \leq j \leq k-1}{\max}\dist_j \\
\sum_{m=k-K}^{k-1}\|c_m\| \nonumber & \leq 2NK\underset{k-K-3\tau \leq j \leq k-1}{\max}\dist_j \\
\sum_{m=k-K}^{k-1}\|d_m\| \nonumber & \leq 2(1+\gamma L)NK\underset{k-K-2\tau \leq j \leq k-1}{\max}\dist_j \\
\sum_{m=k-K}^{k-1}\|Sx_m\|          & \leq 2K\underset{k-K \leq j \leq k-1}{\max}\dist_j\enspace,
\end{align}
the last inequality following from Corollary~\ref{lem::Scocoercive}.

From the definition of $\Sigma_K(k)$ in~\eqref{eq::S(K)} and from~\eqref{eq::boundSumDist}, by introducing
\[
Y:=1+\gamma L+2\beta,
\]
we have that
\begin{align*}
\Sigma_K(k) &\leq \sum_{m=k-K}^{k-1}(\|d_m\| + \beta\|c_m\| + \beta \|b_m\| + \|Sx_m\|) \\
            &\leq \underbrace{2K(YN+1)}_{W(K)}\underset{k-K-3\tau \leq j \leq k-1}{\max}\dist_j
\end{align*}

 Since~\eqref{eq::S(K)} is bounded, we can accordingly bound~\eqref{eq::boundakRecursive}, ~\eqref{eq::boundbkRecursive} and~\eqref{eq::boundckRecursive}:
  \begin{subequations}
  \begin{alignat}{3}
  \|a_k\| \nonumber&\leq \eta N\Sigma_{2\tau}(k) \\
          \nonumber&\leq \eta NW(2\tau) \underset{k-5\tau \leq j \leq k-1}{\max}\dist_j \\
                   &= \eta 4\tau N(YN+1) \underset{k-5\tau \leq j \leq k-1}{\max}\dist_j\label{eq::boundakDist}\\
  \|c_k-b_k\| \nonumber&\leq \eta N2\Sigma_{3\tau}(k) \\
              \nonumber&\leq \eta  N2W(3\tau)\underset{k-6\tau \leq j \leq k-1}{\max}\dist_j \\
                       &= \eta 2(3\tau)N(YN+1) \underset{k-6\tau \leq j \leq k-1}{\max}\dist_j \label{eq::boundbkDist}\\
  \|d_k\| \nonumber&\leq \eta N(1+\gamma L)\Sigma_{2\tau}(k) \\
          \nonumber&\leq \eta N(1+\gamma L)W(2\tau)\underset{k-5\tau \leq j \leq k-1}{\max}\dist_j \\
                   &= \eta(1+\gamma L)4\tau N(YN+1)\underset{k-5\tau \leq j \leq k-1}{\max}\dist_j\enspace. \label{eq::boundckDist}
 \end{alignat}
\end{subequations}

Using~\eqref{eq::boundakDist},~\eqref{eq::boundbkDist} and~\eqref{eq::boundckDist}, $\|e_k\|$ from~\eqref{eq::boundek} can be bounded as
\begin{align}{\label{eq::boundekDist}}
 \|e_k\| \nonumber&\leq \|d_k\| + \beta\|c_k-b_k\| \\
                  &\leq \eta \underbrace{N(YN+1)(4\tau(1+\gamma L)+6\beta\tau)}_{X}\underset{k-6\tau \leq j \leq k-1}{\max}\dist_j\enspace,
\end{align}
where we bounded the quantities with the maximum delay that appeared in~\eqref{eq::boundakDist},~\eqref{eq::boundbkDist} and~\eqref{eq::boundckDist}.\qed

\section{Proof of Thoerem~\ref{thm::inequality_convergence} }{\label{app::Conv5:Convergence}}
{\it Proof}
Note that~\eqref{eq::EquationConvergence} simplifies to
\[
  \eta^2X^2\left(\frac{1}{\epsilon}+\frac{\eta(1+\delta)}{\delta}\right) < \nu-\epsilon \enspace.
\]
As a result, we need to find parameters $\eta,\beta,\gamma,\delta,\epsilon$ such that the following set of inequalities are satisfied:
\begin{equation}{\label{eq::EquationConvergence2}}
\left\{
\begin{array}{l}
 \eta^2X^2\left(\frac{1}{\epsilon}+\frac{\eta(1+\delta)}{\delta}\right) < \nu-\epsilon,\\
 Y = 1+\gamma L+2\beta,\\
 X = N(YN+1)(4\tau(1+\gamma L)+6\beta\tau),\\
 \delta>0,\\
 \epsilon>0,\\
 \beta>0,\\
  0<\gamma<\gamma_{\max},\\
  0<\rho<\frac{1}{2(1+\delta)}.
 \end{array}
 \right.
\end{equation}
The upper bound $\gamma_{\max}$ ensures that the stepsize $\gamma$ is admissible (a possible option is, \eg, $\gamma_{\max}=2/L$ as proven in Appendix~\ref{app::properties}).
We start be noting that the values of $\delta$ and $\epsilon$ are irrelevant as long as they are positive. To this end, we can start by choosing $\epsilon$ such that
$\nu-\epsilon>0$.
From the inequality $\eta<1/(2(1+\delta))$ it follows that
\[
\frac{1}{\epsilon} + \frac{\eta(1+\delta)}{\delta} < \frac{2\delta+\epsilon}{2\delta\epsilon},
\]
thus having
\[
 \eta^2 < \frac{2\delta\epsilon(\nu-\epsilon)}{X^2(2\delta+\epsilon)},
\]
from which the result follows.\qed



\section{Cocoercivity and quasi-strong monotonicity of $S$}{\label{app::properties}}
\subsection*{Proof of Corollary~\ref{lem::Scocoercive}. }{\label{app::Scocoercive}}
{\it Proof}  From~\cite[Proposition~4.33]{book_comb} we have that $T$ is nonexpansive if and only if $S$ is $1/2$-cocoercive. Hence it suffices to show that $T$ is nonexpansive. From Assumption~\ref{Ass::propertiesB}, $B$ is $1/L$-cocoercive, which means that $\gamma B$ is $1/\gamma L$-cocoercive. It follows from~\cite[Lemma~5.1~(iv)]{PrecGFBS} that $T_B=I-\gamma B$ is $\gamma L/2$-averaged. From~\cite[Proposition~4.25~(i)]{book_comb} it follows that $T_B$ is nonexpansive provided that $\gamma<2/L$. Finally, from Assumptions~\ref{Ass::propertiesA} and~\ref{Ass::propertiesB} we conclude that $T$ is nonexpansive as the composition of nonexpansive operators.
\qed

\begin{lemma}\label{lem::Sstrongmonotone}
If $B$ is $\mu$-strongly monotone, then the operator $S$ is quasi-$\nu$-strongly monotone, where $\nu=1-\sqrt{(1-2\gamma\mu+\mu\gamma^2L)}$.
\end{lemma}
{\it Proof}  The Lemma is proven in~\cite[Proposition~2]{peng2015arock} for the case of the proximal gradient method. The proof below is essentially the same generalized for an operator $T$.
 From~\cite[Example~22.5]{book_comb} we have that if $T$ is $c$-Lipschitz continuous for some $c\in[0,1)$ then $I-T$ is $(1-c)$-strongly monotone. Let us then prove that $T$ is indeed Lipschitz continuous. For any $x\in\mathcal{H}$ and $x_\ast\in\operatorname{fix} T$ it holds that:
 \begin{align*}
 \|T_Bx-T_Bx_\ast\|^2 \nonumber&= \|x-x_\ast\|^2-2\gamma\langle x-x_\ast,Bx-Bx_\ast\rangle +\gamma^2\|Bx-Bx_\ast\|^2 \\\nonumber
  &\leq \|x-x_\ast\|^2-\gamma(2-\gamma L)\langle x-x_\ast,Bx-Bx_\ast\rangle \\\nonumber
  &\leq \|x-x_\ast\|^2-\mu\gamma(2-\gamma L)\|x-x_\ast\|^2 \\
  &= (1-2\gamma\mu+\mu\gamma^2L)\|x-x_\ast\|^2\enspace,
\end{align*}
where the first inequality follows from the $1/L$-cocoercivity of $B$, while the second one from the $\mu$-strong monotonicity of $B$.

Thus $\|T_Bx-T_Bx_\ast\|\leq \sqrt{(1-2\gamma\mu+\mu\gamma^2L)}\|x-x_\ast\|$ and since $T_A$ is nonexpansive, we have that $\|Tx-Tx_\ast\|\leq \sqrt{(1-2\gamma\mu+\mu\gamma^2L)}\|x-x_\ast\|$. Finally, $S$ is quasi-$\nu$-strongly monotone with $\nu=1-\sqrt{(1-2\gamma\mu+\mu\gamma^2L)}$ for $\gamma<2/L$.
\qed

\bibliographystyle{spmpsci_unsrt}
\bibliography{biblioAsyncParallelInertialKM_final}

\end{document}